\documentclass[12pt]{amsart}
\usepackage{amscd,amssymb}
\usepackage[english]{babel}
\usepackage[all]{xy}

\newtheorem{Lem}{Lemma \thesection.}
\newtheorem{Th}[Lem]{Theorem \thesection.}
\newtheorem{Cor}[Lem]{Corollary \thesection.}
\newtheorem{Def}[Lem]{Definition \thesection.}

\newtheorem{Prop}[Lem]{Proposition \thesection.}

\begin{document}
\title[CLASS
$\rm {VII_{0}}$ SURFACES]{ CLASS $\rm {VII_{0}}$ SURFACES WITH $b_2$
CURVES
}

\author[DLOUSSKY, OELJEKLAUS \& TOMA]{
Georges  DLOUSSKY, Karl OELJEKLAUS \& Matei TOMA}
\thanks{{\it Georges Dloussky, Karl Oeljeklaus}: LATP-UMR(CNRS) 6632 .
CMI-Universit\'e d'Aix-Marseille I.
39, rue Joliot-Curie, F-13453 Marseille Cedex 13, France. e-mail:
dloussky@gyptis.univ-mrs.fr,  karloelj@gyptis.univ-mrs.fr\\
{\it Matei Toma}: Fachbereich Mathematik-Informatik, Universit\"at Osnabr\"uck,
49069 Osnabr\"uck, Germany and
Institute of Mathematics of the Romanian Academy.
e-mail: Matei.Toma@mathematik.Uni-Osnabrueck.DE}
\date{\today}


\begin{abstract}
{We give an affirmative answer to a conjecture of Ma. Kato, namely that every
compact complex surface $S$  in Kodaira's class $\rm {VII_{0}}$ having
$b_2(S) > 0$ and $b_2(S)$ rational curves, admits a global spherical shell.  }
\end{abstract}

\maketitle

\tableofcontents
A {\it surface} is a compact complex manifold of dimension  2 .
        We denote  by $b_{i}(S)$ the i-th Betti number of  $S$. \newpage
\section{Introduction}
\setcounter{Lem}{0}
Kodaira's class $VII_0$, which consists of minimal compact complex
surfaces $S$ having $b_1(S)=1$, is not completely understood so far. 
In fact, only the case when $b_2(S)=0$
is completely classified, by the work of Kodaira \cite{Kod},
Inoue \cite{In1}, Bogomolov \cite{Bog}, Li-Yau-Zheng \cite{LYZ} and Teleman \cite{T}.\\
For $b_2(S)>0$ a construction method and thus a large subclass of surfaces has been introduced
 by Kato \cite{KA}.
These are exactly the minimal surfaces with $b_2>0$
containing {\it global spherical shells}.
(See the next section for the definition.)
One can show that a surface $S$ of class $VII_0$ has
at most $b_2(S)$ rational curves on it and that if moreover $S$
admits a global spherical shell, then there are exactly
$b_2(S)$ rational curves on $S$.
Ma. Kato conjectured that the converse should be true as well.
Important progress towards this conjecture was made by
I. Nakamura \cite{N1}, \cite{N2}, who showed that if $S$ has $b_2(S)$ rational curves,
then their configuration is that of the curves of a surface with global spherical shells
and $S$ is a deformation of a blown up Hopf surface.\\
This paper is devoted to the proof of Kato's conjecture:\\

{\bf Main Theorem:} If $S$ is a surface of class $VII_0$
with $b_2(S)>0$ and with $b_2(S)$ rational curves,
then $S$ admits global spherical shells.\\

At present all known surfaces $S$ of class $VII_0$ 
with $b_2(S)>0$ contain global spherical shells. In fact,
by making additional assumptions on $S$ like the existence
of a homologically trivial divisor \cite{E} or of two
cycles of rational curves \cite{N1} or of a holomorphic
vector field \cite{DOT1}, \cite{DOT2}, it was shown that
$S$ contains global spherical shells.
The present paper makes a further step in this
direction. \\

The paper is organized as follows.
Section $1$ is preparatory. We recall some facts
on surfaces with global spherical shells and 
on surfaces $S$ of class $VII_0$ with $b_2(S)>0$
rational curves. Such surfaces were called
{\it special} by Nakamura \cite{N2}. We also
prove a fact we shall need later, namely
that if $S$ is special, then the canonical
bundle of a suitable finite ramified covering
of $S$ is numerically divisorial.

Using the knowledge of the configuration of rational curves
on a special surface, we prove in section $2$ the existence of a 
logarithmic $1$-form twisted by a flat line bundle.
Passing to a finite ramified covering we get a global twisted holomorphic vector
field. The twisting is again by some flat line bundle. This induces
a true holomorphic vector field on the universal covering
$\tilde S$ of $S$. \\
In section $3$ we prove that this holomorphic vector field is 
completely integrable and that the universal covering
of the complement of the curves of $S$ is isomorphic to 
$\mathbb H \times \mathbb C$, where $\mathbb H$
denotes the complex half plane.

Section $4$ is devoted to the computation of
the action of the fundamental group on $\mathbb H \times \mathbb C$.

This allows us to recover in section $5$ the contracting rigid
germ of holomorphic map which gave birth to our surface $S$.
Using the work of C. Favre \cite{Fav} which classifies such germs,
we are able to conclude.\\

{\it Acknowledgments:} This paper was written during a stay 
of the third author at the LATP of CNRS, Universit\'e de 
Provence. He wishes to thank this institution for its 
hospitality as well as  the CNRS and the DFG for financial support.

\section{Preliminaries}

We start by recalling some definitions and known facts.
A compact complex surface $S$ is said to belong
Kodaira's class $VII$ if $b_1(S)=1$ and
to class $VII_0$ if it is moreover minimal.
Surfaces of class $VII_0$ with $b_2=0$
have been completely classified, see for instance \cite{T}.
In this paper we deal with the case $b_2>0$.
It is well known then, that the Kodaira dimension
of $S$ is negative and that the algebraic dimension
vanishes. In particular, $S$ has finitely many
irreducible curves.\\
At present the only known surfaces of class $VII_0$
with $b_2>0$ contain {\it global spherical shells (GSS)}.
A global spherical shell is a neighborhood $V$ of
$S^3 \subset \mathbb C^2 \setminus \{0\}$ which is
holomorphically embedded in the surface $S$ such that
$S \setminus V$ is connected. All surfaces with GSS
may be constructed by a procedure due to Ma. Kato \cite{KA}.
As a consequence they all have exactly $b_2(S)$ rational
curves; some of them admit an elliptic curve as well.
 Ma. Kato also made the following\\

{\bf Conjecture} If $S$ is a class $VII_0$ surface with $b_2>0$ rational curves,
then $S$ admits a GSS.\\
 
Following I. Nakamura \cite{N2} we shall call a 
class $VII_0$ surface {\it special}, if it has $b_2>0$ rational curves.
Since special surfaces admitting homologically trivial
divisors or with all rational curves organized in one or two
cycles have been shown to admit GSS, \cite{E}, \cite{N1}, \cite{D2}, \cite{KA}, we concentrate
our attention on the remaining ones. We shall call them
{\it special surfaces of intermediate type}.\\

In \cite{N2}, Nakamura proved that the configuration of the rational
curves of a special surface of intermediate type is the same as
that of a surface with GSS with the same $b_2$. In particular the dual graph of such a configuration 
is connected and contains a cycle to which some trees are attached. 
Nakamura also showed that these surfaces are deformations
of blown-up primary Hopf surfaces, in particular that their fundamental
group is isomorphic to $\mathbb Z$. Furthermore, he proposed
a line of attack to Kato's conjecture. However, one of his 
conjectures proved to be not correct, see \cite{Zaf}.\\

{\it Notation} We denote by $D$ the maximal reduced divisor
of a special surface $S$, by $M(S)$ the intersection matrix of the rational curves of $S$ and set 
$k(S):=\sqrt{|\det M(S)|}+1$.
(It is well known that $\sqrt{|\det M(S)|}$  is
the index of the subgroup
generated by  classes of  curves in $H^2(S, \mathbb Z)$ and thus it is an integer).
 Moreover we denote by $\tilde S$ the universal cover of $S$ and by $\tilde D$ 
the preimage of $D$ in $\tilde S$. Then $\tilde D$ is the universal cover of $D$.\\

In this section we show that after passing to a ramified
covering and a resolution of singularities, we may suppose
that the anticanonical bundle of our surface is numerically
divisorial, i.e. there exists a flat line bundle $L$
on $S$ and a divisor $D_{-K}$ such that 
$K_S^{-1} \simeq L \otimes \mathcal O(D_{-K})$.

Since under our assumptions $\pi_1(S) \simeq \mathbb Z$,
it is easy to see that the flat line bundles are parametrized 
by $\mathbb C^* \simeq Hom(\pi_1(S),\mathbb C^*) \simeq Pic^0(S)$.
We shall often write $L^{\lambda}$ for the line  bundle
corresponding to the complex number $\lambda \in \mathbb C^*$.

\begin{Lem} Let $S$ be a special surface of intermediate type.
Then there exist a positive integer $m$, a flat line bundle $L$
and an effective divisor $D_m$ such that
$$( K_S \otimes L)^{\otimes m} = \mathcal O(-D_m).$$
\end{Lem}
{\it Proof}: Since the cohomology classes associated to the
rational curves of $S$ generate $H^2(S,\mathbb Q)$, there always
exist $m \in \mathbb N^*$, $L \in Pic^0(S)$ and a divisor
$D_m$ such that
$$( K_S \otimes L)^{\otimes m} = \mathcal O(-D_m).$$
We have only to check that $D_m \geq 0$.
Let $D_m= D_+ - D_-$ with $D_+,D_- \geq 0$ and $D_+.D_- \geq 0$.
The adjunction formula implies that $D_m.C \leq 0$ for any
irreducible curve $C$ on $S$. Hence
$$0 \geq D_m.D_- = (D_+ - D_-).D_- \geq -D_-^2 \geq 0.$$
But then $D_-=0$, since $S$ does not admit homologically trivial divisors.
\qed

\begin{Def}{\rm The smallest possible $m \in \mathbb N^*$
for which a decomposition  
$$( K_S \otimes L)^{\otimes m} = \mathcal O(-D_m)$$
as in Lemma \thesection.1 exists, will be called
 the {\it index} of the surface $S$ and denoted by $m(S)$. When $m(S)=1$ we denote $D_{-K}=D_1$ and call this
 the {\it numerically anticanonical divisor} of $S$.}
\end{Def}
Notice that $D_m$ is unique when $S$ is of intermediate type. 
The following proposition will enable us to reduce the
proof of the Main Theorem to the case of special
surfaces of index $1$. We have formulated it for simplicity for special surfaces of intermediate type, but 
the general case can be proved similarly.

\begin{Prop} \label{RevRam} Let $S$ be a special surface of intermediate type  
with index $m:=m(S)>1$.  Then there exists a diagram
$$
\xymatrix{
 & S' \ar[dr]^{\pi'}   &  \\
T \ar[ur]^c  \ar[d]_{\rho} & & Z'\\
Z \ar[dr]_{\pi}& & T' \ar[u]_{\rho'} \ar[dl]^{c'}\\
&S& 
}
$$
where
\begin{itemize}
\item $(Z,\pi,S)$ is 
a m-fold cyclic ramified covering space  of $S$, branched over $D_m$, 
\item  $(T,\rho,Z)$ is
 the minimal desingularization of $Z$,
\item $(T,c,S')$ is the contraction of the (possible) exceptional 
curves of the first kind,
\item $S'$ is a special surface 
with action of the group of $m$-th roots of unity
$\mathbb U_m$, with index $ m(S')=1$,
\item $Z'$ is the quotient space of $S'$ by $\mathbb U_m$, 
\item $(T',\rho',Z')$ is the minimal desingularization of $Z'$,
\item $(T',c',S)$ is the contraction of the (possible) 
exceptional curves of the first kind,
\end{itemize}
such that the restriction over $S\setminus D$ is commutative, i.e.
$$\theta:=\pi\circ \rho\circ c^{-1}
=c' \circ {\rho'}^{-1} \circ \pi' : S'\setminus D' \to S\setminus D.$$
Moreover, if $S'$ has a GSS, then $S$ has a GSS as well.
\end{Prop}

Proof: We have
$$( K_S \otimes L)^{\otimes m} = \mathcal O(-D_m).$$

Let $X$ be the total space of the line bundle $K^{-1} \otimes L^{-1}  $. 
We choose an open trivialisation covering 
$\mathcal U=(U_i)$ for $K$ and $L$
with local coordinates $(z_1^i,z_2^i)$, defining cocycles $(k_{ij})$ and $(g_{ij})$ of 
$K$ and $L$ such that 
$D_m\cap U_i=\{f_i=0\}$ and 
$$k_{ij}^m g_{ij}^m\frac{f_i}{f_j}=1.$$
If $\zeta_i$ is the fiber variable of $X$ over $U_i$, 
the equations $\zeta_i^ m=f_i(z)$ fit together and define an
analytic subspace $Z\subset X$.   
It is easy to see that $\mathbb U_m$ 
acts holomorphically and effectively on $Z$ and that $Z/\mathbb U_m=S$. 
Let  $\pi:Z\to S$ be the  projection on $S$.  
The ramified covering $(Z,\pi,S)$ is  branched 
exactly over $supp(D_m)$. 
The local meromorphic 
2-forms 
$$\omega_i=\frac{dz_1^i\wedge dz_2^i}{\zeta_i}$$
yield a twisted meromorphic 2-form $\omega$ on $X$, 
hence on $Z$, for 
$$\omega_i=\frac{dz_1^i\wedge dz_2^i}{\zeta_i}
=\frac{k_{ij}^{-1}dz_1^j\wedge dz_2^j}{(k_{ij} g_{ij})^{-1}\zeta_j}=
g_{ij}\frac{dz_1^j\wedge dz_2^j}{\zeta_j}=g_{ij}\omega_j.$$
Now, let $(T,\rho,Z)$ be the minimal desingularization of $Z$. This includes of course the
normalization of $Z$. Notice that the normalization is connected by the minimality of $m$. 
Set $H=\rho^\star L$, 
then  $\tau=\rho^\star\omega $ is a twisted meromorphic 2-form on $T$ 
which does not vanish and has a non trivial polar divisor $E$ and $K_T \otimes H= \mathcal O(-E)$. 
Finally, in order 
to obtain $S'$, we contract all exceptional curves of the first 
kind. It is clear that the index $S'$ is one. As before $\mathbb U_m$ 
acts holomorphically on $S'$
and $S'\setminus D'$ 
is a covering manifold of $S\setminus D$.  

The quotient $Z'=S'/\mathbb U_m$ is a normal surface. 
The desingularization $(T',\rho',Z')$ yields $S$ 
after contraction of the exceptional curves of the first kind.\\
Since $S$ has no non-constant meromorphic functions, 
the same holds for $S'$. Applying then the 
classification of Kodaira, $S'$ is a  
K3 surface, a torus or a surface of class VII$_0$. The 
first two cases have a trivial canonical bundle, 
hence $S'$ belongs to class VII$_0$. Moreover $b_2(S')>0$, 
because it contains a cohomologically non trivial divisor.\\

We shall now show that if $S$ has a GSS, then $S'$ will
also have this property. We may choose a GSS 
$S^3 \subset V \subset \mathbb C^2 \setminus \{0\}$
which after embedding in $S$ cuts only one curve $C$ of our surface $S$.
For a suitable choice of global coordinates $(z_1,z_2)$
on $V$, the intersection $C \cap V$ is given by the equation $z_1=0$.
Moreover $K_S|_V$ and $L|_V$ are trivial on $V$. Let $\zeta$ be the fiber coordinate of $X$ over $V$. 
Suppose $D_m \cap V = n(C \cap V)$. 
The pull-back $\pi^{-1}(S^3)$ of the sphere $S^3$ to $Z$ is given by the equations
$$|z_1|^2 +|z_2|^2 =1, \ \ z_1^n=\zeta^m.$$
Let $d:= g.c.d.(n,m)$ and $n'=n/d,\ m' = m/d$.
Then there are $d$ irreducible components around $\pi^{-1}(S^3)$ in $Z$
which will become disjoint after normalization.
We denote by $\Sigma_1,...,\Sigma_d$ the corresponding components
of $\pi^{-1}(S^3)$. Let us choose $\Sigma_1$ with equations
$$|z_1|^2 +|z_2|^2 =1, \ \ z_1^{n'}=\zeta^{m'}.$$ 
This component is normalized by the map 
$$(t,z_2) \mapsto (t^{m'},z_2,t^{n'}).$$
The pull-back $\Sigma'_1$ of $\Sigma_1$ to $T$ will be given in these 
coordinates by the equation 
$$|t|^{2m'} + |z_2|^2 =1.$$
We check now that $T \setminus \Sigma'_1$ is connected.
Let $P,Q \in T \setminus \Sigma'_1$ two points in a neighborhood
of $\Sigma'_1 \cap \rho^{-1}(\pi^{-1}(D))$ which find themselves on 
different sides of $\Sigma'_1$. Their projections $\pi(\rho(P))$, $ \pi(\rho(Q))$
on $S$ may be connected by a path avoiding $S^3$ and $D$. We may
lift this path to a path $\gamma$ in $Z$ which connects
$\rho(P)$ to some $Q'\in \pi^{-1}(\pi(\rho(Q)))$. In order to ensure
that $Q'$ and $\rho(Q)$ coincide, it is enough to let the initial
path in $S$ turn the needed number of times around the components
of $D$. By further lifting $\gamma$ to $T$ we get the looked for
connectedness. For $m'>1$ ,
$$ \Sigma'_1 = \{(t,z) \in \mathbb C^2 \ \ |\ \ |t|^{2m'} + |z_2|^2 =1\}$$
is no longer a sphere, but remains the border of a bounded Stein domain
of $\mathbb C^2$. Noticing that $\Sigma'_1$ may be approximated by a strictly pseudoconvex hypersurface (for example of equation $\epsilon|t|^2+ |t|^{2m'} + |z_2|^2 =1$) it is possible to use the same arguments as in 
\cite{D1} in order to get a contracting holomorphic germ and hence a GSS
on $S'$.\\

In particular, when $S$ has a GSS, $S'$ has exactly $b_2(S')= -K_{S'}^2$
rational curves. But since for a special surface $S$, the dual graph of the curves
is the same as for some surface with a GSS and since the intersection matrix of
the curves of $S'$ depends only on this graph, we get that our $S'$
also has $b_2(S')= -K_{S'}^2$ rational curves. Thus $S'$ is special
under the weaker condition that $S$ is special (of intermediate type).

Finally, it is not difficult to show (cf. \cite{D3}),
that quotients of GSS surfaces by the actions of finite cyclic groups
of automorphisms remain GSS surfaces. Thus if $S'$ has a GSS,
then $S$ will also have one.  \qed

\section{Existence of a twisted logarithmic $1$-form} 
\label{s1form}
\setcounter{Lem}{0}
In this section we shall consider a special surface $S$ 
 of  intermediate type and prove that it always admits a twisted logarithmic 1-form. 
 Let $D=D_{max}= \sum_{i=1}^n C_i$ be the maximal
reduced divisor of $S$. 

\begin{Lem}\label{L1}
If $S$ has index $m(S)=1$, 
then 
the numerically anticanonical divisor satisfies $D_{-K} > D$. In particular, for every flat line bundle 
$L$ on $S$, $H^0(S,K_S\otimes L)=0$.
\end{Lem}

{\it Proof}: 
Suppose  that $D_{-K}$ does not contain all the curves of $S$.
Since by \cite{N2} the maximal divisor $D$ is connected, there will exist
an irreducible curve $C$ which is not contained in $D_{-K}$ but such that
${\rm supp}(D_{-K}) \cap C \neq \emptyset$. Then $K_S.C  = -D_{-K}.C <0$
which gives a contradiction to the adjunction formula. Thus $D_{-K} \geq D$. 
But the equality $D_{-K} = D$
would imply $D^2 = - b_2(S)$ which may happen only on 
Inoue-Hirzebruch surfaces by \cite{N2}. But Inoue-Hirzebruch surfaces
are not of intermediate type, so 
in our case we  have 
$D_{-K}>D$. \qed

\begin{Lem}\label{L2}
For every $L \in {\rm Pic}^0(S)$ we have 
$\Gamma(S,\Omega^1 \otimes L)=0$.
\end{Lem}

{\it Proof} : We may suppose $L\neq 0$.
Using the exact sequence of sheaves
$$ 0 \rightarrow d{\mathcal O}(L) \rightarrow \Omega^1 \otimes L \rightarrow 
\Omega^2 \otimes L \rightarrow 0$$

 we get $$\Gamma(\Omega^1 \otimes L)=
\Gamma(d{\mathcal O}(L)).$$

Take now $\omega \in \Gamma(d{\mathcal O}(L))$ and denote
by $\tilde \omega$ its pull-back on the universal covering $\tilde S$ of $S$.
Then $g^*\tilde \omega = \lambda \tilde \omega$ where $g$ is a generator
of $\pi_1(S)\simeq \mathbb Z$ and $\lambda$ is the twisting factor which
corresponds to $L \in {\rm Pic}^0(S) \simeq \mathbb C^*$.
Let $f$ be a primitive of $\tilde \omega$ and $c \in \mathbb C$ such that 
$f \circ g = \lambda f +c$. Replacing $f$ by $h:= f + \frac{c}{\lambda -1}$,
we get $h \circ g = \lambda h$ which means that $h$ induces a section in 
$\Gamma(S,\mathcal O(L))$. By our assumption on $S$, $h$ has then to be the
zero section and thus $\omega = 0$. \qed

\begin{Lem}\label{L3}
A non-trivial twisted logarithmic $1$-form on $S$ is always closed 
and has poles along each curve of $S$.
\end{Lem}

{\it Proof} : Let $0 \neq \omega \in \Gamma(S,\Omega^1({\rm log}D) \otimes L)$
for a flat line bundle $L$. Then $d\omega 
\in \Gamma(S,\Omega^2(D) \otimes L).$ If $d\omega \neq 0$,
then its associated divisor $\Gamma$ satisfies $0 \leq -\Gamma \leq D$ which 
contradicts Lemma \thesection.\ref{L1}. 
 Thus $d\omega =0$.
By Lemma \thesection.\ref{L2}, the pole divisor $D_{\infty}$
of $\omega$ is non-trivial. Then $D_{\infty}$ must contain the cycle of rational
curves of $S$, otherwise one could write $\omega$ as a non-trivial logarithmic
$1$-form in a neighborhood $V$ of $D_{\infty}$ and since the dual graph of
 $D_{\infty}$ would be  contractible, $\omega$ would be holomorphic on $V$
by \cite{SS}. So let now $C_1$ be an irreducible component of $D_{\infty}$
and suppose there exists a rational curve $C_2$ not contained in $D_{\infty}$
such that $\emptyset \neq C_1 \cap C_2 = \{p\}$.
Choose coordinates $(z_1,z_2)$ locally around $p$ such that $C_i = \{z_i = 0\}$,
$i=1,2$, and write $$\omega = \alpha_1 \frac{dz_1}{z_1} 
+\alpha_2 dz_2$$ in these coordinates. Since $$0=d\omega =
\frac{\partial\alpha_1}{\partial z_2} dz_2 \wedge \frac{dz_1}{z_1}
+\frac{\partial\alpha_2}{\partial z_1} dz_1 \wedge dz_2,$$ 
we see that $\alpha_1$ must have the form 
$$\alpha_1(z_1,z_2)=\beta(z_1)+ z_1 z_2 \gamma(z_1,z_2).$$
If $\beta(0)\neq 0$, the restriction of $\omega$ to $C_2$ would have as only pole
a simple pole in $p$, which is impossible.
Therefore $\beta(0)= 0$ and $$\omega = (\frac{\beta}{z_1} +z_2 \gamma)dz_1
+\alpha_2 dz_2$$ has no pole along $C_1$, which gives a contradiction.
Thus $D_{\infty} = D$. \qed

\begin{Lem}\label{L4}
Two logarithmic $1$-forms $\omega_1,\omega_2$ on $S$ twisted by flat
line bundles $L_1, L_2$ are necessarily linearly dependent.
\end{Lem}
{\it Proof} : If $\omega_1,\omega_2$ were not linearly dependent,
their exterior product $\omega_1 \wedge \omega_2 \in \Gamma(S, \Omega^2(D)
\otimes L_1 \otimes L_2)$ would be non-identically zero.
This contradicts Lemma \thesection.\ref{L1}. \qed

\begin{Lem}\label{L5}
If $L$ is a flat line bundle on $S$ and $\Gamma(S, \Omega^1({\rm log}D)
\otimes L^{-1}) = 0$, then the morphism 
$$ H^2(S,\mathbb C_S(L)) \rightarrow H^2(D,\mathbb C_D(L))$$
induced by restriction is bijective.
\end{Lem}
{\it Proof} : We may suppose $L$ to be non-trivial, since otherwise
the conclusion holds by our hypotheses on $S$. 

 
The long exact cohomology sequence of the diagram 

$$
\xymatrix{
0 \ar[r] & \mathbb C_D(L) \ar[r] & \mathcal O_D(L) 
\ar[r]^{d \ \ \ \  } & \bigoplus_{i=1}^n 
\Omega_{C_i} \otimes L \ar[r] & 0 \\
0 \ar[r] & \mathbb C_S(L) \ar[u]_{\rm restr} \ar[r] & 
\mathcal O_S(L) \ar[u]_{\rm restr} 
\ar[r]^{d} & d
\mathcal O_S(L) \ar[u]_{\rm restr} \ar[r] & 0 & (1) 
}
$$

gives 

$$
\xymatrix{
H^1(\mathcal O_D(L)) \ar[r] & H^1( \bigoplus_{i=1}^n \Omega_{C_i} \otimes L)
\ar[r] & H^2(\mathbb C_D(L)) \ar[r] & H^2(\mathcal O_D(L)) \\
 \ar[u] H^1(\mathcal O_S(L)) \ar[r] & 
H^1(d \mathcal O_S(L)) \ar[u] 
\ar[r] & H^2(\mathbb C_S(L))  \ar[u] \ar[r] & H^2(\mathcal O_S(L)) \ar[u].  
}
$$

Using the identities
$$H^2(\mathcal O_S(L))= 0,\quad H^0(\mathcal O_S(L))=0,\quad H^2(\mathcal O_D(L))=0,
\quad H^0(\mathcal O_D(L))=0$$
 and the theorem of Riemann-Roch, we get 
$H^1(\mathcal O_S(L))=0$ and $H^1(\mathcal O_D(L))=0$. Thus our task comes
to showing the bijectivity of the morphism
$$H^1(d \mathcal O_S(L)) 
\rightarrow H^1(\bigoplus_{i=1}^n \Omega_{C_i} \otimes L).$$
On the other hand the commutative triangle	
$$
\xymatrix{
&\bigoplus_{i=1}^n 
\Omega_{C_i} \otimes L & \\
d \mathcal O_S(L) \ar[ur]^{\rm restr} \ar[rr]& & \Omega^1(L) \ar[ul]_{\rm restr},
}
$$
the long exact cohomology sequence associated to

$$
  0 \rightarrow d \mathcal O_S(L)\rightarrow \Omega^1_S(L) \rightarrow
\Omega^2_S(L) \rightarrow 0   \leqno{(2)}     
$$

and the fact that $H^0(\Omega^2_S(L))=0$, allow us to get the commutative triangle

$$
\xymatrix{
&H^1(\bigoplus_{i=1}^n 
\Omega_{C_i} \otimes L) & \\
H^1(d \mathcal O_S(L)) \ar[ur] \ar[rr]& & H^1(\Omega^1(L)) 
\ar[ul],
}
$$
where the horizontal arrow is an isomorphism. Thus, we have only to 
prove that the morphism

$$  H^1(\Omega^1(L)) \rightarrow H^1(\bigoplus_{i=1}^n 
\Omega_{C_i} \otimes L)   \leqno{(3)}$$
is bijective.
We examine the dimensions first.
$$ \dim H^1(\bigoplus_{i=1}^n \Omega_{C_i} \otimes L) =
\sum_{i=1}^n \dim H^1( \Omega_{C_i} \otimes L) =n$$
and
$$\dim H^1(\Omega^1(L)) = - \chi(\Omega^1(L)) + h^0(\Omega^1(L)) 
+h^2(\Omega^1(L))=  - \chi(\Omega^1(L)) = - \chi(\Omega^1) =b_2,$$
since we know that $h^0(\Omega^1(L))=h^0(\Omega^1(L^{-1})))=0$.
By assumption $b_2 = n$, so it is enough to show only
the surjectivity of the morphism (3). In order to do so, let us
compute the kernel $\mathcal N$ of the surjective morphism
$$ \Omega_S^1 \rightarrow \bigoplus_{i=1}^n \Omega_{C_i}.$$

Locally around a point $C_1 \cap C_2$, where $C_i=\{z_i =0\}$,
this morphism is given by 
$$f_1dz_1 + f_2dz_2 \mapsto (f_1(z_1,0)dz_1,f_2(0,z_2)dz_2).$$
Thus a section of the kernel must have the form $z_2g_1dz_1+z_1g_2dz_2$
where $g_1,g_2$ are holomorphic functions.
Therefore we get the duality

$$
\xymatrix{
{\mathcal {N}} \otimes \Omega^1_S(log D) \ar[r]^{\ \ \ \ \ \ \  \wedge} &  
\Omega^2_S }
$$

$$
\xymatrix{
(z_2 g_1 dz_1 +z_1 g_2 dz_2,\displaystyle \frac{h_1}{z_1} dz_1+\frac{h_2}{z_2} dz_2) 
\ar@{|-{>}}[r] 
& (g_1 h_2 - g_2 h_1) dz_1 \wedge dz_2 , 
}
$$
proving that 
 
$$\mathcal N \simeq \Omega^1_S(log D)^{\vee} \otimes \Omega^2_S. $$

In order to finish the proof we consider now the long exact
cohomology sequence of

$$ 0 \rightarrow \mathcal N \otimes L \rightarrow \Omega^1_S \otimes L
\rightarrow \bigoplus_{i=1}^n \Omega_{C_i} \otimes L \rightarrow 0$$

and use the fact that $H^2(\mathcal N \otimes L) = H^0(\Omega^1_S(log D)
\otimes L^{-1})^* = 0 $ which is ensured by  hypothesis. $\qed$

\begin{Th}\label{FORME}
Let $S$ be a special surface of intermediate type and $k:=k(S)$.
Then there is a choice of a generator of $\pi_1(S) \simeq \mathbb Z$
such that $S$ admits a closed logarithmic $1$-form twisted by the flat
line bundle $L^k$.
\end{Th}

{\it Proof} : By \cite{N2} there exists a surface $S'$ with a GSS
such that the dual graph of the maximal reduced divisor $D'$ of $S'$
coincides with that of $D$. By \cite{DO} and \cite{Fav}, Thm. 1.2.24,
we may choose the generator
of $\pi_1(S') \simeq \pi_1(D')$ such that $\Gamma(S',\Omega^1(log D')
\otimes L^k) \neq 0$. We further fix the generator of $\pi_1(S)\simeq 
\pi_1(D) \simeq \mathbb Z \simeq \pi_1(S') \simeq \pi_1(D')$ to be the same
as above. We may now suppose that $\Gamma(S,\Omega^1(log D)
\otimes L^{1/k}) = 0$, otherwise we change the generator of $\pi_1(S)$.
We shall identify sections of sheaves on $S$ twisted by $L^k$ with sections
of the pullback-sheaves on the universal cover $\tilde S$ of $S$ which
respect the representation $\rho:\pi_1(S) \rightarrow \mathbb C$
defining $L^k$. We start with the exact sequence 
$$ \ \ 0 \rightarrow d \mathcal O(L^k) \rightarrow 
d \mathcal O(log D) \otimes L^k \rightarrow \bigoplus_{i=1}^n 
\mathcal O_{\tilde C_i}
\otimes L^k \rightarrow 0 \leqno{(4)} $$
where $\tilde C_i$ are the normalisations of the curves $C_i \subset S$
and with a certain element 
$$ a \in \Gamma(S,\bigoplus_{i=1}^n \mathcal O_{\tilde C_i} \otimes L^k)=
\Gamma(\tilde S,\bigoplus_{i=- \infty}^{\infty} \mathcal O_{\tilde C_i})^{\rho}$$
to be defined below. Seen as an element in $$\Gamma(\tilde S,\bigoplus_{i=- \infty}^{\infty} 
\mathcal O_{\tilde C_i}) \simeq 
\bigoplus_{i=- \infty}^{\infty} \Gamma(\tilde S, \mathcal O_{\tilde C_i}), $$
$a$ becomes a vector in $\mathbb C^{\mathbb Z}$ with the property 
$ a_{i+n} = k a_i$ for all $i \in \mathbb Z$.\\
Choose a non-trivial element

$$ \omega' \in \Gamma(S',\Omega^1(logD') \otimes L^k)$$
and put
$$a_i:= \int_{\gamma_i}  \omega',$$
where $\gamma_i$ is a small path around $C_i' \subset \tilde D'$ and we denote again by $  \omega'$  the pull-back of $\omega'$ to the universal cover
 $\tilde S'$ of $S'$. 
The Camacho-Sad formula for the foliation defined by $ \tilde \omega'$
gives
$$ C_j^2 = - \sum_{{i \neq j} \atop {C_i \cap C_j \neq \emptyset}} 
\frac{a_i}{a_j}, $$
see \cite{DOT2}.
Moreover we have seen in \cite{DOT2} that $\omega'$ may be chosen such that
$a \in \mathbb Z[1/k]^{\mathbb Z}$. For such a choice let further $U_j$
be small neighborhoods of the curves $C_j$ on $\tilde S$ and consider
divisors $D_j$ in these neighborhoods of the form
$$D_j := k^{\nu}(a_jC_j + \sum_{{i \neq j} \atop {C_i \cap C_j \neq \emptyset}}
a_i C_i),$$
for some $ \nu \in \mathbb N$ which is sufficiently large.
Since $D_j.C_j = 0$, $D_j$ is the zero divisor of some holomorphic function $f_j
\in \mathcal O(U_j)$.
Put then $$ \omega_j :=k^{-\nu} \frac{df_j}{f_j} \in \Gamma(U_j, 
d \mathcal O(log \tilde D)).$$
One may choose the functions $f_j
\in \mathcal O(U_j)$ such that 
$$g^*\omega_{j+n} = k \omega_j$$
holds for all $j \in \mathbb Z$; here $g \in \pi_1(S)$ denotes 
the ``positive'' generator of $ \pi_1(S)$. 
There exist local coordinate functions $z_j,z_i$ 
 which define the curves $C_j$, $C_i$
such that
around $C_i \cap C_j $ we have 
$$ \omega_j = a_j \frac{dz_j}{z_j} + 
a_i \frac{dz_i}{z_i} $$
and around points of $C_j$ away from any other compact curve
 $$ \omega_j = a_j \frac{dz_j}{z_j} .$$
 We look now at the situation on
$U_i \cap U_j$. Here $\omega_j-\omega_i$ is a closed holomorphic $1$-form
and thus exact if $U_i \cap U_j$ is simply connected, which we shall
always suppose; in fact, $\omega_j-\omega_i$ has the form 
$$ a_j \frac{dg_j}{g_j} +a_i \frac{dg_i}{g_i} ,$$
with $g_i,g_j \in \mathcal O^*(U_i \cap U_j)$. \\
Let $\mathcal U$ be the covering of $\tilde S$ which consists
of the open sets $U_i$ and of $\tilde S \setminus \tilde D$.
We set $\omega =0$ on $\tilde S \setminus \tilde D$. Then we get a 
cocycle $(\omega_j-\omega_i)_{j,i} \in H^1(\mathcal U, d \mathcal O(L^k))$
which represents the image of $a$ through the canonical connecting homomorphism
associated to the sequence (4). We follow this image further through the
isomorphism 
$$ H^1(d \mathcal O(L^k)) 
\overset{\simeq}{\longrightarrow} 
H^2(\mathbb C(L^k))$$
which comes from the short exact sequence
$$0 \rightarrow \mathbb C(L^k) \rightarrow \mathcal O(L^k)\rightarrow 
d \mathcal O(L^k) \rightarrow 0.$$    
In order to do this we pass to a finer covering $\mathcal V = (V_{\nu})_{\nu}$
which has the property that the sets $V_{\mu \nu} := V_{\mu} \cap V_{ \nu}$
are simply connected. Let $\gamma $ be the refinement map and set 
$\omega_{\nu}:= \omega_{\gamma(\nu)}|_{V_{\nu}}$ for 
$ V_{\nu} \subset U_{\gamma(\nu)}$. 
Since $\omega_{\nu} -\omega_{\mu}$ is exact
on $V_{\mu \nu}$, we may choose primitives $f_{\nu \mu}$. This is done such that
$$ f_{\nu \mu}=0 \ \ \  {\rm if}  \ \ \  \gamma(\mu)= \gamma(\nu)\; {\rm and} $$
$$ f_{\nu \mu}=f_{\gamma(\nu) \gamma( \mu)} 
\ \ \  {\rm if}  \ \ \  \gamma(\mu)\neq \gamma(\nu)\; {\rm and}\; \gamma(\mu),\, \gamma(\nu) \in \mathbb Z.$$
The cocycle 
$$((f_{\mu \nu} - f_{\lambda \nu} + f_{\lambda \mu})
|_{V_{\lambda \mu \nu}})_{\lambda \mu \nu}$$
is the looked for  element in $H^2(\mathcal V, \mathbb C(L^k))$.
But we remark now that the trace of this cocycle on $\tilde D$
is zero and thus its image in $H^2(D \cap \mathcal V, \mathbb C(L^k))$
by the restriction morphism will be zero as well. One sees this
by checking the different possibilities for $\lambda, \mu, \nu,$ 
such that $V_{\lambda \mu \nu} \cap \tilde D \neq \emptyset :$
\begin{itemize}
\item when $\gamma(\lambda)= \gamma(\mu)= \gamma(\nu)$, all primitives
are zero
\item when $\gamma(\lambda)= \gamma(\mu) \neq  \gamma(\nu)$, one gets again
$f_{\mu \nu} - f_{\lambda \nu} + f_{\lambda \mu}|_{V_{\lambda \mu \nu}}=0$.
\end{itemize}
Notice that at least two elements among 
$\gamma(\lambda), \gamma(\mu), \gamma(\nu)$ must be equal.

 Now Lemma \thesection.\ref{L5} shows that the class of our cocycle
in $H^2(S,\mathbb C(L^k))$ vanishes. Hence there exists a non-trivial
element $\omega$ in $\Gamma(S, d \mathcal O(log D) \otimes L^k)$ which maps
onto $a \in \Gamma(S, \bigoplus_{i=1}^n \mathcal O_{\tilde C_i} \otimes L^k)$.
\qed

\begin{Cor}\label{C1}
The foliation associated to $\omega \in \Gamma(S, \Omega^1(log D) \otimes L^k)$ 
has simple singularities and rational Camacho-Sad indices with respect
to the curves of $S$.
\end{Cor}
{\it Proof} : If $\omega$ is obtained as in the main part of
the proof of Theorem \thesection.\ref{FORME} then the associated Camacho-Sad indices are 
$\frac{a_i}{a_j}$, i.e.
the same as those associated to $\omega'$ on $S'$, and thus rational.
But if $\omega \in \Gamma(S, \Omega^1(log D) \otimes L^{1/k})$, where
the orientation on $\pi_1(S)$ and $\pi_1(S')$ is chosen to be the same, we 
have to reconsider the Camacho-Sad equations of the components $C_i$
of $D$. Let $\mathcal F$ denote the foliation associated to $\omega$.
Lemma \thesection.\ref{L3} implies that the irreducible components of $D$ are
invariant for $\mathcal F$. Let $C_1, C_2$ be two such components which
intersect at $p$ and $(z_1,z_2)$ local coordinate functions around $p$
such that $C_i = \{z_i = 0\}$ for $i=1,2$. 
We may write $\omega$ around $p$ as
$$ \omega = g_1\frac{dz_1}{z_1}+g_2\frac{dz_2}{z_2}$$
where $g_1,g_2$ are holomorphic functions in $z_1,z_2$.
Consider now small paths $\gamma_i$ turning around
$C_i$, contained say in the local curve $z_{3-i}=c_{3-i}$,
for two constants $c_1,c_2 \in \mathbb C$.
Since $\omega$ is closed the integrals

$$
\int_{\gamma_1}\omega =  2 \pi i Res_{z_1 =0}(\omega|_{z_2 = c_2})
= 2 \pi i g_1(0,c_2) $$
$$\int_{\gamma_2}\omega = 2 \pi i Res_{z_2 =0}(\omega|_{z_1 = c_1})
= 2 \pi i g_2(c_1,0) 
$$
are independent of $c_1,c_2$.
Moreover, since $\omega$ has true poles of order one along $C_1$ and $C_2$, both
integrals are non-zero; in particular $g_1(0,0) \neq 0 \neq g_2(0,0)$
and $p$ is a simple singularity for $\mathcal F$. Now the foliation
$\mathcal F$ is defined locally around  $p$ by the kernel of the
form $z_2 g_1 dz_1 + z_1 g_2 dz_2$ and the Camacho-Sad index
of $\mathcal F$ with respect to $C_2$ is by definition
$$ CS(\mathcal F,C_2,p):= Res_{z_1=0}\bigl( \frac{\partial}{\partial z_2}
(- \frac{z_2 g_1}{z_1 g_2})|_{C_2}\bigr)
=Res_{z_1=0}\bigl( -\frac{g_1(z_1,0)}{z_1 g_2(z_1,0)}\bigr)= $$ $$
= -\frac{g_1(0,0)}{g_2(0,0)}= - \int_{\gamma_2}\omega/\int_{\gamma_1}\omega $$
and similarly for $C_1$.
A consequence of this is that all Camacho-Sad indices of
$\mathcal F$ along the curves of $S$ are non-zero.\\

Now recall that the maximal divisor $D$ of $S$ consists of a cycle
$\sum_{i=1}^m C_i$ of rational curves to which a non zero
number of trees of rational curves is attached. It is easy
to see that the Camacho-Sad indices associated to the curves
of the trees are all rational. Indeed, if $B_{i,1},...,B_{i,m_i}$
is the tree with root $C_i$, and if $- B_{i,j}^2=:b_{i,j}$, we get:
$$ CS(\mathcal F, B_{i,m_i}, B_{i,m_i} \cap B_{i,m_i -1}) =- b_{i,m_i}$$
for the top,
$$ CS(\mathcal F, B_{i,m_i -1}, B_{i,m_i} \cap B_{i,m_i -1}) =-\frac{1}{b_{i,m_i}},$$
$$ CS(\mathcal F, B_{i,m_i-1}, B_{i,m_i -1} \cap B_{i,m_i -2}) =- b_{i,m_i -1}
+\frac{1}{b_{i,m_i}}$$
and so on.
Setting as in \cite{DOT2} 
$b_i : = -{\rm CS}({\mathcal F},C_i,C_i \cap B_{i,1})$,
$ d_{i}:= -C_{i}^{2} - b_{i}$ and 
$ \alpha_{i}:= -{\rm CS}({\mathcal F},C_i,C_{i-1} \cap C_i), \ \
i=1,\ldots,m $, we get the equations  
$$
\alpha_{i}+\displaystyle \frac{1}{\alpha_{i+1}} = d_{i}, \ \ {\textrm {for }} 
{1 \leq i \leq m-1}
$$
and
$$
\begin{array}{rcl}
\alpha_{m}+\displaystyle \frac{1}{\alpha_{1}} = d_{m}.\\
\end{array}
$$
Each $\alpha_{i}$ is the solution of a quadratic equation with
rational coefficients and since we know already that a rational solution
exists by working with $\omega'$ on $S'$ , the other solution
has to be rational as well. \qed

\begin{Cor}
A special surface of intermediate type of index 1
 possesses a non-trivial holomorphic vector field $\theta$
twisted by some flat line bundle and its vanishing divisor
$D_{\theta}$ contains all the curves of  $D$ except perhaps
the summits of the trees. In particular $\theta$ vanishes on all the curves
of the cycle of $D$.
  
\end{Cor}
{\it Proof} : Let $S$ be a special surface of intermediate type
and let $D_{-K}$ be the (unique) numerical anticanonical divisor
on $S$. This means that there exists a torsion factor $\kappa \in
\mathbb C^*$ such that
$$ \mathcal O(D_{-K}) = { K}_S^{-1} \otimes L^{\kappa}.$$
Let now $0 \neq \omega \in \Gamma(S,\Omega^1(log D)
\otimes L^{k})$ and $Z$ the subspace of the intersection points
of the curves of $S$. By Corollary \thesection.\ref{C1} this
is exactly the space of singularities of the associated
foliation to $\omega$. Thus $\omega$  induces an exact sequence
$$ 0 \rightarrow \mathcal O(-D) \otimes L^{1/k}\rightarrow
\Omega^1 \rightarrow L^{\kappa k} \otimes {\mathcal J}_Z\otimes\mathcal O(D-D_{-K})
\rightarrow 0 $$
and by duality we get a non-trivial section 
$$\theta \in \Gamma(S, \Theta_S \otimes L^{k \kappa} )$$
vanishing on $D_{\theta} := D_{-K} -D$.
We know by Lemma \thesection.\ref{L1} that $D_{\theta} > 0$.
Let now $C_1 \not\subset Supp(D_{\theta})$ such that $C_1$
intersects an irreducible curve $C_2 \subset D_{\theta}$.
Since $\theta$  defines the same foliation as $\omega$
the curve $C_1$ is $\theta $- invariant. On the other hand
since $C_1 \cap C_2$ is a singularity of this foliation and $\theta$
vanishes on $C_2$ the vanishing order of the restriction of $\theta$
to $C_1$ is at least two. But then since $C_1$ is rational this order
is exactly two and $C_1$ cannot intersect another curve of $D$.
Thus $C_1$ is the top of a tree of rational curves.
This implies our statement. \qed

\section{The universal covering of the complement of the curves}
\setcounter{Lem}{0}
From now on we shall consider a special surface $S$ of intermediate
type which admits a numerical anticanonical divisor. We have
seen that in this case $S$ possesses a non-trivial twisted
logarithmic $1$-form $\omega$ with twisting factor $k=k(S)$
and a non-trivial twisted holomorphic vector field $\theta$
with twisting factor, say $\lambda \in \mathbb C^*$.
The case $\lambda = 1$, i.e. $\theta$ is a holomorphic
vector field on $S$, was considered and completely understood
in \cite{DOT1} and \cite{DOT2}. In this section we  prove  that 
in the general case, as in case 
$\lambda = 1$, the
universal covering of $S\setminus D$ is isomorphic to
$\mathbb H \times \mathbb C$, where $\mathbb H$ 
denotes the complex half-plane.

Let $U$ be a small open neighborhood of $D$, such that
$D$ is a deformation retract of $U$.
We have $$\pi_1(U)=\pi_1(D)=\pi_1(S)=\mathbb Z$$
and we denote as before by $g$  a generator of this group. 
There is a fundamental domain $U_0$ for the action of
$\mathbb Z$ in the inverse image $\tilde U$
of $U$ in the universal cover  $\tilde S$ of $S$,
 such that the border of $U_0$
in $\tilde U$ cuts $\tilde D$ in a component $C_0$
and in its translated $g(C_0)$ along a circle $S^1$.
Set $Y_0:= \bigcup_{\nu \geq 0} g^{\nu}(U_0)$.
We keep the notation $\omega$ for the logarithmic 1-form one gets on $\tilde S$. 

\begin{Lem}

There is a normalization of $\omega$ such that the representation
$$ \rho : \pi_1(\tilde S \setminus \tilde D) \rightarrow \mathbb C$$

 $$\gamma \mapsto \int_{\gamma} \omega$$

\noindent has as image $2 \pi i \mathbb Z[\frac{1}{k}] \subset \mathbb C$.
Furthermore, one can choose this normalization such that 
$\rho(\pi_1(Y_0 \setminus \tilde D)) = 2 \pi i \mathbb Z$.

\end{Lem}

{\it Proof} : Since $\tilde S$ is simply
connected, the group $\pi_1(\tilde S \setminus \tilde D)$ is
generated by small paths around the irreducible components 
of $\tilde D$. 
Thus keeping the notations of the previous section, 
we see that $\rho(\pi_1(\tilde S \setminus \tilde D))$ is a 
$\mathbb Z[\frac{1}{k}]$-module generated by 
$$2\pi i a_0= \int_{\gamma_0} \omega ,...,2\pi i a_{n-1}= \int_{\gamma_{n-1}} \omega,$$
where $\gamma_0,...,\gamma_{n-1}$ are small paths around the curves
$C_0,...,C_{n-1}$ in $U_0$.

Using  Corollary \ref{s1form}.\ref{C1} and the way we computed the
Camacho-Sad indices, we see that  we can normalize the form $\omega$ such that the numbers
$a_0,...,a_{n-1}$ are non-zero integers 
with g.c.d. equal to $1$;
so $\rho(\pi_1(\tilde S \setminus \tilde D))$ is free of rank $1$
as $\mathbb Z[\frac{1}{k}]$-module.
 
In a similar way the group 
$\rho(\pi_1( Y_0 \setminus \tilde D))$
is generated as a $\mathbb Z$-module 
by small paths around the irreducibles components 
of $\tilde D$ which meet $U_0$ and we get the announced result. \qed\\

In what follows we suppose $\omega$ to be normalized such that
$\rho(\pi_1(Y_0 \setminus \tilde D))=2 \pi i \mathbb Z$.
This is the same as to say that $a_0,...,a_{n-1}$ are
non-zero integers with g.c.d. equal to $1$.

Let $A$ be a fundamental domain for the action of $\mathbb Z$
on $\tilde S$ and $X_0:= \bigcup_{j \geq 0} g^j( A)$. 
Translating by $g$, we may suppose that
$\tilde D \cap X_0 \subset Y_0$. We remark that after such a translation
we have
$$\rho(\pi_1((Y_0 \cup X_0) \setminus \tilde D))=
\rho(\pi_1(Y_0 \setminus \tilde D)) = 2 \pi i \mathbb Z.$$

Fix $z_0 \in U_0$. We define a holomorphic function $f$ on
$(Y_0 \cup X_0) \setminus \tilde D$ by
$$ f(z) = exp \bigl( \int_{z_0}^z \omega + \frac{1}{k-1} 
\int_{z_0}^{g(z_0)} \omega \bigr).$$

One verifies easily that $f$ is well defined and that
$$f(g(z))= f^k (z)  \leqno{(\dag)}$$
for $z \in (Y_0 \cup X_0) \setminus \tilde D$.\\ 

Let $C$ be the smooth part in $\tilde D$ of an irreducible component 
of $\tilde D \cap Y_0$. Since $\omega$ is a closed logarithmic
$1$-form, an easy computation shows that one can extend $f$
meromorphically across $C$ such that $C$ belongs to the
zero or to the polar set of the extension of $f$. By replacing 
$\omega$ with $-\omega$ if necessary, one finds at least one
preimage on $\tilde S$ of a component of the cycle of $D$
in the zero-set of $f$. This means that the connected components of the 
polar set of $f$ are exeptional divisors in $Y_0$. Thus the polar set
is empty, see \cite{SS}, and $f$ vanishes on $\tilde D \cap Y_0$.   In particular, we see that the integers 
$a_0,...,a_{n-1}$ are positive.  
\begin{Lem}\label{LL1}
$|f(z)| < 1 $ for any $z \in Y_0 \cap X_0$.
\end{Lem}
{\it Proof} :          Remark first that one can extend the function
$|f|$ to the whole of $\tilde S$. The extended function
is still denoted by $|f|$. Then the statement of the Lemma
can be rephrased by saying that the image of $|f|$ is
the interval $[0,1[$. Suppose now that this is not so.
Then $|f|$ has as image $[0,\infty[$ since $|f| \circ g= |f|^k$. We consider the real
hypersurface $\tilde H := |f|^{-1}(1)$, which is $\pi_1(S)$-invariant
and thus descends to a compact real hypersurface $H$ on $S$.
Obviously, $H$ is a (compact) leaf of the foliation defined
by $\Re e \, \omega$. Remark that the morphism 
$\pi_1(H) \rightarrow \pi_1(S)$ is non-trivial. Otherwise
the connected components of $\tilde H$ would be compact
and their intersection with $f^{-1}(1)$ would give compact analytic curves
in the complement of $\tilde D$, which is absurd. By passing to a
finite unramified covering of $S$, we may even suppose that
$\pi_1(H) \rightarrow \pi_1(S)$ is surjective.\\
Next we  prove that $\omega$ must have non-vanishing
periods on $\tilde H$. If not, consider a $\pi_1(S)$-invariant neighborhood
$\tilde V$ of $\tilde H$ which is the preimage of a neighborhood
$V$ of $H$ in $S$ and on which $\omega$ has no periods;
thus $\omega |_{\tilde V}$ is exact. We define
the following holomorphic function on $\tilde V$:
$$h(z):= \int_{z_0}^z \omega + \frac{1}{k-1} \int_{z_0}^{g(z_0)} 
\omega, $$
where $z_0 \in \tilde V$ is a fixed base point, $z \in \tilde V$ and 
integration is done along paths in $\tilde V$. We have $h(g(z)) = k h(z)$
for all $z \in \tilde V$.
The function $h$ does not take the value $0$, because the set
$\{h=0\}$ would then be a $\pi_1(S)$-invariant analytic curve giving
rise to a compact curve in the complement of $D$ in $S$, which
is absurd. Hence one may consider the $1$-form
$d(log |h|)$ which is closed and descends to a closed non-twisted
form on $V$. This form obviously defines the same foliation
as $\Re e \omega$. But this means that this foliation has trivial
holonomy. This implies that the leaves near $H$ are also compact.
In particular they are contained in $V$. The inverse images 
of these leaves are leaves of the foliation defined by $|f|$ on $\tilde S$.
On the other hand they are completely contained in $\tilde V$ and intersect
all sets $g^{\nu}(\tilde S \setminus X_0)$, $\nu \in \mathbb Z$.
But the relation $f(g(z)) = f^k(z)$ and the choice of $V$ imply
that
$$ \bigcap_{\nu \in \mathbb Z} |f|(g^{\nu}(\tilde V \setminus X_0)) =\{1\},$$
giving a contradiction.\\
Thus $\omega|_{\tilde H}$ has non-vanishing periods. The next point is to show
that $f$ may be extended holomorphically to $$U_1 := |f|^{-1}(]1, \infty[).$$
In order to do this we consider the sets 
$$U_{\alpha} := |f|^{-1}(]\alpha, \infty[)$$ for $\alpha \geq 1$ and
$$ M:= \{\alpha > 1 | f {\rm \ admits\  a\  holomorphic\  extension\  to\ } 
U_{\alpha}\}.$$
Then $M$ is non-empty, since $U_{\alpha} \subset X_0$ for $\alpha > 
\underset{z \in \partial X_0}{ \max} |f(z)|$. 
Furthermore,  the set $M$ is closed in $]1, \infty[$ since 
$U_{\alpha_0} = \bigcup_{\alpha >\alpha_0} U_ {\alpha}$ for each ${\alpha_0}$.
It remains to check that $M$ is open as well. Let $\alpha \in M$. There exists some $\nu \in 
\mathbb N $ such that $ U_{\alpha} \subset g^{- \nu}(X_0)$. On $g^{- \nu}(X_0)$
one can extend $f^{k^{\nu}}$ holomorphically. Take now a finite open covering
$V_1,...,V_{\mu}$ of $(\partial U_{\alpha}) \setminus X_0$ such that each $V_i$
intersects $U_{\alpha}$ and that on $V_i$ a $k^{\nu}$-th root of 
$f^{k^{\nu}}$ is defined. Consider next on each $V_i$ that $k^{\nu}$-th root
of $f^{k^{\nu}}$ which coincides with $f$ on $V_i \cap U_{\alpha}$.
This gives an extension of $f$ to 
$$U_{\alpha} \cup X_0 \cup \bigcup_{i=1}^{\mu} V_i.$$
Remark that this set will contain some $U_{\beta} $ with $\beta < \alpha$.\\
We finish the proof of the Lemma by considering a path  $\gamma \subset U_1$
such that $\int_{\gamma} \omega \neq 0$, which is possible since
$\omega$ has non-vanishing periods on $H$. But for $\nu$ sufficiently large 
$$\int_{g^{-\nu} \circ \gamma} \omega 
= \frac{\int_{ \gamma}\omega}{k^{\nu}},$$
 cannot be a multiple of $2\pi i$, which is incompatible with the definition of $f$
on the whole of $U_1$. \qed

\begin{Prop} The holomorphic vector field $\tilde \theta \in \Gamma(\tilde S, \tilde
\Theta_{\tilde S})$ induced by $\theta$ is completely integrable on
$\tilde S$.
\end{Prop}
{\it Proof} : We consider the integrability of $\tilde \theta$ along a leaf
$F$ of the foliation. It will be enough to show that one can find a local
integration radius which is uniform for all $F$. We set $A_0:=
\overline{X_0 \cup Y_0} \setminus g(X_0 \cup Y_0)$. The set $A_0$ intersects
$\tilde D$  along a curve $C_0$ and its translated by $g$, $C_n$.
There is a constant $c \in ]0,1[$ such that $F \subset |f|^{-1}(c)$.
Using the relation $|f| \circ g= |f|^k$ again and Lemma \thesection.\ref{LL1}
we see that $\lim_{\nu\to\infty}|f|\circ g^\nu=0$ uniformly on compact sets, hence there exists a $\nu_0 \in \mathbb N$ such that
$F \cap \bigcup_{\nu > \nu_0} g^{\nu}(A_0) = \emptyset$.
On the other side there exists a $\nu_1 \in \mathbb N$ such that 
$F \cap \bigcup_{\nu \geq \nu_1} g^{-\nu}(A_0) \subset \tilde U$.
By passing to a suitable translation we may assume that $\nu_1 =0$.
We may further assume that around $\tilde U \cap \partial A_0$ 
we have coordinate functions $(z_0,z_1)$ such that $C_0 = \{z_0=0\}$,
$$\omega = a_0 \frac{dz_0}{z_0}$$ with  $a_0 \in \mathbb N^*$ and 
$$\tilde \theta = \alpha(z_0,z_1)z_0^s \frac{\partial}{\partial z_1}$$
with $s \in \mathbb N^*$ and $\alpha$ a nowhere 
vanishing holomorphic function.
Here we have $f(z_0,z_1) = z_0^{a_0}$.  Thus $$(z_0,z_1) \in
F \cap \partial(\bigcup_{\nu \geq 0} g^{\nu}(A_0))$$
implies $|z_0|^{a_0} = c $. But since 
$F \cap \bigcup_{\nu \geq 0} g^{\nu}(A_0)$ is contained in the
compact set $\bigcup_{\nu = 0}^{\nu_0} g^{\nu}(A_0)$, the integration radius
of $\tilde \theta$ at points of 
$$F \cap \bigcup_{\nu \geq 0} g^{\nu}(A_0)$$ is at least as large as
the integration radius of $\tilde \theta$ at 
$$F \cap  \partial(\bigcup_{\nu = 0}^{\nu_0} g^{\nu}(A_0))
=F \cap  \partial(\bigcup_{\nu \geq 0} g^{\nu}(A_0))$$ and this is
the minimal integration radius of $\alpha(z_0,z_1) z_0^s \frac{\partial}
{\partial z_1}$ at points $(z_0,z_1)$ with $|z_0|^{a_0} = c $.
Looking now at $\tilde \theta$ on 
$$F \cap  \partial(\bigcup_{\nu \geq -r } g^{\nu}(A_0))$$
for $ r \in \mathbb N$, we see by applying $g_*^r$, that the integration
radius at these points will be at least as large as the minimal integration radius
of $ \lambda^r \alpha(z_0,z_1) z_0^s \frac{\partial}
{\partial z_1}$ at points $(z_0,z_1)\in
F \cap \partial(\bigcup_{\nu \geq 0} g^{\nu}(A_0))$ 
with $|z_0|^{a_0} = c^{k^r} $.
But the sequence of velocities 
$$ |\lambda|^r  c^{{sk^r}/a_0} \sup |\alpha|$$ is obviously bounded 
and thus there is an uniform integration radius for $\tilde \theta$
on $F$. \qed\\

The kernel $\ker \rho$ defines a covering 
$$ \pi : X' \rightarrow \tilde S \setminus
\tilde D.$$ One checks immediately that the $\mathbb Z$-action induced
by $g$ on $\pi_1(\tilde S \setminus \tilde D)$ stabilizes $\ker \rho$
and thus induces a $\mathbb Z$-action on $X'$. We denote again
by $g$ a lift of $g$ on $X'$. Thus we get an action of the
semi-direct product $\mathbb Z \ltimes \mathbb Z[\frac{1}{k}]$
on $X'$ whose quotient is $S \setminus D$. Let $\omega' = \pi^*
\omega$ and $\phi: X' \rightarrow \mathbb C$ be a primitive
of $\omega'$ on $X'$, such that $\exp (\phi)$ and $f \circ \pi$
coincide on a connected component of the $\pi$-preimage of
$(Y_0 \cup X_0) \setminus \tilde D$. 
Since $\phi \circ g =k \phi$, the image of $\phi$ is invariant under
the action of the multiplicative group $\{k^{\nu}| \nu \in \mathbb Z\}$;
this image is also invariant under the action of the additive
group $2 \pi i \mathbb Z[\frac{1}{k}]$. Since $f$ takes its values in the
unit disk by Lemma \thesection.\ref{LL1}, we see now that $\phi(X')$
must coincide with the left half plane $\mathbb H_l:= \{w \in \mathbb C |
\Re e(w) <0 \}$. The function $\phi : X' \rightarrow \mathbb H_l$
is a holomorphic surjection. The connected components of its fibers 
are leaves of the inverse image foliation induced by $\mathcal F$.
For the proof of the next Proposition the reader is referred to \cite{DOT2},
Prop. 2.2. Note that the Camacho-Sad indices of our foliation here,
are positive integers, just as in \cite{DOT2}. 
\begin{Prop}
The fibers of $\phi$ are connected.
\end{Prop}
\begin{Cor}
The foliation defined by $\theta$ on $\tilde S \setminus \tilde D$
has no closed leaf.
\end{Cor} 
{\it Proof} : If $F$ were a closed leaf on $\tilde S \setminus \tilde D$,
its preimage $\pi^{-1}(F)$ in $X'$ would also be closed. But since
$\rho|_{\pi_1(F)}$ is trivial, the group $\mathbb Z[\frac{1}{k}]$
operates on $X'$ by permuting components of $\pi^{-1}(F)$. Using
the preceeding Proposition, the non-discreteness of the
$2 \pi i\mathbb Z[\frac{1}{k}]$-orbits in $\mathbb H_l$ and the fact that $\phi$ is a submersion, we get
a contradiction. \qed
\begin{Lem} The fibers of $\phi$ are isomorphic to $\mathbb C$.
\end{Lem}
{\it Proof} : Since $\tilde \theta $ is completely integrable,
there is a holomorphic $\mathbb C$-action on $\tilde S \setminus \tilde D$.
The fixed point set of a non-trivial element of $\mathbb C$ is a closed
analytic subset of $\tilde S \setminus \tilde D$ which is a union of
leaves of the foliation. By the preceeding Corollary such a fixed point set
cannot have dimension $1$. Hence it is either empty or the whole space
$\tilde S \setminus \tilde D$. We must exclude the second case.
In this situation the $\mathbb C$-action factorizes to a $\mathbb C^*$-action.
Thus all fibers of $\phi$ are isomorphic to $\mathbb C^*$. Moreover
the $\mathbb C^*$-action lifts to $X'$ making $\phi : X' \rightarrow \mathbb H_l$
into a principal $\mathbb C^*$-bundle. But such a bundle over $\mathbb H_l$
is always trivial. We may therefore see $\phi$ as the first factor projection
$$ X' \simeq \mathbb H_l \times \mathbb C^* \rightarrow \mathbb H_l.$$
If $(w,z)$ denote coordinate functions on $\mathbb H_l \times \mathbb C^*$
the pull-back of the vector field has the form $\alpha(w)z \frac{\partial}{\partial
z}$ for some $\alpha \in \mathcal O^*(\mathbb H_l)$. We now consider
generators $f_g,f_{\gamma}$ of the groups $\mathbb Z$ and $\mathbb Z[\frac{1}{k}]$
acting on $\mathbb H_l \times \mathbb C^*$.

By passing to a double covering of $S$ we may assume that  $f_g$
 acts on the $\mathbb C^*$-fibers by  homotheties. Suppose that the same holds for $f_{\gamma}$.
Then we have
$$ f_g(w,z) =(kw, \beta(w)z)$$
$$ f_{\gamma}(w,z)=(w + 2 \pi i, \gamma(w)z)$$
for some $\beta,\gamma \in \mathcal O^*(\mathbb H_l)$.
The compatibility with the pulled-back vector field
implies
$$ \alpha(kw)=\lambda \alpha(w),\, 
\alpha(w + 2 \pi i)=\alpha(w)$$
for all $w \in \mathbb H_l$. 
The second relation implies that $\alpha$ is the composition of a function
$u$ on the punctured unit disc with the exponential. Then the first relation translates into
$$ u(\zeta^k)=\lambda u(\zeta ), \, {\rm for \, all}\, \zeta \in {\Delta}^{\star};$$
(by $\Delta$ we denoted the unit disk in $\mathbb C$).
But this implies $\lambda=1$. Hence an effective $\mathbb C^*$-action  
on $S$, which is excluded by \cite{Hau}.

When $f_{\gamma}$ is composed with an inversion a similar argument applies,
working with $\alpha^2$ instead of $\alpha$ for instance.

Thus the $\mathbb C$-action is effective and the fibers of $\phi$ are 
isomorphic to $\mathbb C$. \qed\\

Using now the lift of the $\mathbb C$-action on $X'$, we get 
 a $\mathbb C$-principal bundle structure on $X'$ over
$\mathbb H_l$. Again, such a bundle is holomorphically trivial.
In conclusion we have proven the following

\begin{Th}The universal covering of $S \setminus D$ is isomorphic
to $\mathbb H_l \times \mathbb C$.
\end{Th}

\section{The action of the fundamental group}
\setcounter{Lem}{0}
We consider a system of  holomorphic coordinates
$(w,z)$ in $\mathbb H_{l} \times \mathbb C \simeq \widetilde {S \setminus D}$.
The integrable vector field 
induced here by $\theta$ has no zeros and is therefore constant
on each fiber of the projection 
of $\mathbb H_{l} \times \mathbb C$
on $\mathbb H_{l}$.
Consequently, this vector field is of the form
$\alpha(w) \frac{\partial}{\partial z}$ on $\mathbb
H_{l} \times \mathbb C$, where $\alpha \in {\mathcal O}^{\star}(\mathbb
H_{l})$. Conjugating by the automorphism 
$$(w,z) \mapsto (w,
\alpha^{-1} \cdot z),$$ one gets $\alpha \equiv 1$.

Let $\gamma \in \pi_{1}(\tilde S \setminus \tilde D)$ be a closed path in
$\tilde S \setminus \tilde D$ with
$\rho(\gamma) = 2\pi i $. \\
Denoting by $g_{\gamma}$ 
the automorphism of $\mathbb H_{l} \times \mathbb C \simeq \widetilde {S \setminus D}$
corresponding to $\gamma$,
we have:

\begin{eqnarray*}
g_{\gamma}(w,z)&=&(w+ 2\pi i, z + f_{\gamma}(w))\\
g(w,z)&=& (kw, \lambda z+ f_{g}(w)).
\end{eqnarray*}

The case $\lambda =1$ was treated in \cite{DOT2}. From now on we shall
therefore assume that $\lambda \in \mathbb C \setminus \{ 0, 1\}$.\\
 
The automorphism $g_{\gamma}$ generates an action of $\mathbb Z$
on $\mathbb H_{l} \times \mathbb C$, which induces a holomorphic $\mathbb
C$-principal bundle

$$\mathbb H_{l} \times \mathbb C / \mathbb Z \rightarrow  \mathbb
H_{l} / \mathbb Z \simeq  {\Delta}^{\star}.$$
 The triviality of this bundle proves the existence
of a holomorphic function  $h: \mathbb H_{l}\rightarrow \mathbb C$
such that

$$h(w + 2 \pi i) - h(w) = f_{\gamma}(w)$$
and conjugation by $(w,z) \mapsto (w, z + h(w))$ gives us
the new form

$$g_{\gamma}(w,z) = (w + 2 \pi i, z).$$

In what follows we suppose therefore that $f_{\gamma} \equiv 0$.
 
We have

$$ g \circ g_{\gamma} \circ g^{-1} = g_{\gamma}^k,$$

which gives a group isomorphism 

$$ <g_{\gamma}, g> \simeq \mathbb Z 
\ltimes \mathbb Z[\frac{1}{k}] \simeq \pi_1(S \setminus D 
) $$

on the one hand, and the $2\pi i-$periodicity of the function $f_{g}$
on the other hand.\\

Factorizing by $\exp : \mathbb H_{l}  \rightarrow
{\Delta}^{\star},\ \  w \mapsto e^w=: \zeta ,$
gives a Laurent series expansion     

$$f_{g}(w) = \sum_{m \in \mathbb Z} a_{m} e^{mw} =
\sum_{m \in \mathbb Z} a_{m} {\zeta}^{m}.$$

A conjugation by

$$ (w,z) \mapsto (w, z + \beta(w))$$

where $\beta$ is a $2 \pi i$-periodic function on $\mathbb
H_{l}$, does not change
the form of $g_{\gamma}$, but replaces $f_{g}$ by

$$ w \mapsto f_{g}(w) + \beta(kw) - \lambda \beta(w) .$$

Let $$h(\zeta):= \sum_{m \in \mathbb Z} a_{m} \zeta^m, \ \
h_{+}(\zeta):= \sum_{m >0} a_{m} \zeta^m. $$

The series $\sum_{l=0}^{\infty} \lambda^{-(l+1)} h_{+}(\zeta^{k^l})$ converges
uniformly on compact sets in ${\Delta}^{\star}$.
To see this it is enough to write $h_{+}(\zeta) = \zeta(\zeta^{-1}
h_{+}(\zeta))$ and to remark that $\zeta^{-1} h_{+}(\zeta)$ is
holomorphic in $0$. Let 
$$\chi(\zeta):= \sum_{l=0}^{\infty} \lambda^{-(l+1)}  h_{+}(\zeta^{k^l}).$$
We have 
$$\lambda \chi(\zeta)-\chi(\zeta^k)=h_{+}(\zeta).$$

One sets $\displaystyle \beta(w) := \chi(e^w)+ \frac{a_0}{\lambda -1} $ and 
gets

$$f_{g}(w) + \beta(kw) -\lambda \beta(w) = \sum_{m<0} a_{m}e^{mw}.$$

We can therefore suppose that $f_{g}(w) = h(e^w)$, where
$h \in \mathcal O(\mathbb P_1(\mathbb C) \setminus \{0\})$ and $h(\infty)=0$.
One still has the possibility to conjugate with $(w,z) \mapsto (w,z +
\chi(e^w))$, where 
$\chi \in \mathcal O(\mathbb P_1(\mathbb C) \setminus \{0\})$. \\

For a
function $$h(\zeta) = \sum_{m<0}a_{m} \zeta^m$$ and $l
\in \mathbb N^{\star}$, one considers $$h_{l}(\zeta) := \sum_{k^l \vert
m, \ m<0} a_{m}\zeta^m $$ and $$h_{0} := h.$$

Each $h_{l}$ is of the form $h_{l}(\zeta)= f_{l}(\zeta^{k^l})$,
for a certain function $f_{l}$.
Let $\chi$ be the formal series $-\sum_{l \geq 1}\lambda^{l-1} f_{l}$.
Formally one has:

$$\begin{array}{lcl}
h(\zeta) + \chi(\zeta^k) -\lambda \chi(\zeta)&=& 
(h(\zeta) - f_{1}(\zeta^k)) + 
\lambda ( f_{1}(\zeta)- f_{2}
(\zeta^k)) +\lambda^2 (f_{2}(\zeta)- f_{3}(\zeta^k))+ \ldots \\
&=& \sum_{l \geq 0}\lambda^{l-1}  (f_{l}(\zeta)- f_{l+1}(\zeta^k))
\end{array}$$

\noindent and each term contains in its  Laurent series expansion
in $\zeta$ only terms $b_{m}\zeta^m$ with
$k \nmid m$.

For $0<R<1,\ l\geq 1 $ and $\vert \zeta \vert \geq R$ one has :

$$ \vert f_{l}(\zeta) \vert \leq \sum_{k^l \mid m \atop { m<0}} 
\vert a_{m} \vert
R^{\frac{m}{k^l}} \leq \sum_{k^l \mid m \atop { m<0}} \vert a_{m} \vert
R^{m + k^l -1} \leq R^{k^l -1} \sum_{m \leq 0} \vert
a_{m}\vert R^m, $$

which shows that the series defining $\chi$
is uniformly convergent on compact sets of 
$\mathbb P_1(\mathbb C) \setminus \{0\}$.
One gets the following normal form for $h$:

$$ h(\zeta) = \sum_{\scriptstyle m < 0
                                   \atop \scriptstyle k \nmid m}
                                   a_{m}\zeta^m.$$

Remark that if $h$ is of this form, any modification
of $h(\zeta)$ by conjugation gives

$$ h(\zeta) + \chi(\zeta^k) - \lambda \chi(\zeta)$$

\noindent which is in normal form  
if and only if the function
$\zeta \mapsto \chi(\zeta^k) -\lambda \chi(\zeta)$
is identically zero.
To see it, write the power series expansions of these 
functions:
$$ \chi(\zeta) = \sum_{m< 0} b_{m} \zeta^m\ , \ - \chi(\zeta^k)+
\lambda \chi(\zeta) = \sum_{m< 0}c_{m} \zeta^m.$$

If there is $c_{r} \neq 0$ for  $r \in \mathbb Z_{-}$, then $k \nmid r$,
$\lambda b_{r}= c_{r}, \ \ \lambda b_{kr}= c_{kr} + b_{r} =  
\lambda^{-1} c_{r}$,
$\lambda^3 b_{k^{2}r}=\lambda^2 c_{k^{2}r} +\lambda^2 b_{kr} = c_{r}, \ \ \ldots$
and the series $\sum_{m\leq 0}b_{m} \zeta^m$ would not be
convergent on $\mathbb P_1(\mathbb C) \setminus \{0\}$, which one verifies
by putting $\zeta = \vert \lambda \vert$ or $\zeta =1$. 
A contradiction!

We are now in the following
normalized situation:

\begin{Prop}
The action of the fundamental group $\pi_1(S \setminus D) \simeq \mathbb Z \ltimes
\mathbb Z[1/k]$ on the universal cover $\mathbb H_l \times \mathbb C$
of $S \setminus D$ is generated by the two automorphisms:
\begin{eqnarray*}
g_{\gamma}(w,z)&= &(w + 2 \pi i, z)\\
g(w,z) &=& (kw,\lambda z + f_{g}(w))
\end{eqnarray*}
\noindent with $f_{g}(w) = h \circ exp(w) =
 H \circ exp(-w)$ where 
$$H(\zeta) = \sum_{\scriptstyle m > 0 \atop 
\scriptstyle k \nmid m} A_{m}\zeta^m,$$
$A_m=a_{-m}$.
\end{Prop}


The elements $g_{l,n}:= g^{-n} \circ g_{\gamma}^l \circ g^n$
for $ n \in \mathbb N, \ l \in \mathbb Z$ form a  subgroup $\Gamma$
of $ \pi_{1}(S \setminus D)$ isomorphic to $\mathbb Z[\frac{1}{k}]$.
Explicitly we have
$$ g_{l,n}(w,z)= (w+\frac{2 \pi i l}{k^n}, z + \sum_{j=0}^{n-1}\lambda^{-j-1}
(f_{g}(k^j w)-f_{g}(k^j w + \frac{2 \pi i l}{k^{n-j}}))).$$

with $f_{g}(w)= H \circ \exp(-w)$ where $H$
is of the form $$H(\zeta)= \sum_{\scriptstyle m > 0
                                   \atop \scriptstyle k \nmid m}
                                   A_{m}\zeta^m.$$
We know that $\Gamma$ acts properly discontinuously. 
Therefore $H$ is non-trivial. 
One verifies easily that if $n\geq m$, 
$g_{p,n}\circ g_{q,m}=g_{p+qk^{n-m},n}.$\\

For $l=1$, let
$$G_{n,j}(\zeta)= H(\zeta^{k^j})-
H(\zeta^{k^j}\exp({-2\pi ik^{j-n}})), \quad 0\leq j<n,$$
$$F_n(\zeta)= \sum_{0\leq j<n} \lambda^{-j-1} G_{n,j}(\zeta).$$

With these notations and putting $\zeta = \exp(-w)$, we have

$$\sum_{j=0}^{n-1} \lambda^{-j-1} (f_{g}(k^j w)-f_{g}(k^j w 
+ \frac{2 \pi i}{k^{n-j}}))= 
\sum_{j=0}^{n-1}\lambda^{-j-1}G_{n,j}(\zeta)=F_n(\zeta),$$
and 
$$g_{1,n}(w,z)= (w+\frac{2 \pi i}{k^n}, z +F_n(\zeta)).$$

The rest of this section is devoted to the proof of  the following 

\begin{Th}
\label{bori}
The function $H$ is a non-constant polynomial.
\end{Th}

This theorem generalizes the similar statement of \cite{DOT2} whose proof
we owe A. Borichev.\\
In what follows we use the notations:\\
$r\mathbb T$ (resp. $r \mathbb D$) the circle (resp. the open disc)
of radius $r>0$ and $K$ the compact set $\{z\in\mathbb C \mid 
3\leq \mid z\mid\leq 3^k\}$. The different circles $r\mathbb T$ are equipped
with the normalized Lebesgue measure $dm(\zeta)$, for which
$\int_{r\mathbb T} dm(\zeta)=1$. \\
If $f$ is holomorphic on $3 \bar{\mathbb D}$, then $wind(f)$ denotes
the number of zeros of $f$ in $3 \mathbb D$.

The proofs of the following lemmas are left to the reader; 
see however \cite{DOT2} for the case $\lambda =1$.

\begin{Lem}\label{Obs1} Let $f$ be a  holomorphic function on the closed
unit disc 
$3^k\bar{\mathbb D}$ such that for all $z\in K$
 $\mid f(z)\mid \geq 
1$. Then one has the relations
$$\int_{3{\mathbb T}}\ln \mid f(\zeta )\zeta ^{-wind(f)}\mid 
dm(\zeta ) = \int_{3^k{\mathbb T}}\ln \mid f(\zeta )\zeta ^{-wind(f)}\mid 
dm(\zeta ),$$
$$\int_{3{\mathbb T}}\ln \mid f(\zeta )\mid 
dm(\zeta ) = \int_{3^k{\mathbb T}}\ln \mid f(\zeta )\mid 
dm(\zeta ) - (k-1) wind(f) \ln 3.$$
\end{Lem}

\begin{Lem}\label{Obs2} For all $a, b\in{\mathbb C}$, 
such that $\mid a\mid >2$ 
and $\mid b\mid <1$, one has
$$\ln \mid a+b\mid \geq \ln(\mid a\mid) -\mid b\mid .$$
\end{Lem}

\begin{Lem} \label{Obs3} Let $f$ be a holomorphic function  
defined in a neighborhood 
of the closed disc $3\bar{\mathbb D}$. We set
$$A=\int_{3{\mathbb T}}\ln ^+ \mid f(\zeta )\mid dm(\zeta ).$$
Hence there is a constant $C>0$, independent of $f$, for which 
$$\ln \int_{2{\mathbb T}}\mid f(\zeta )\mid dm(\zeta ) \leq CA.$$
(Here $\ln^+ : = \max(\ln,0)$.)
\end{Lem}

\begin{Lem}\label{Obs4} Let $f=\sum_{s\geq 0}\hat f(s)z^s$ be a 
holomorphic function 
in a neighborhood of the disc $2\mathbb D$. 
Then, for all $s\in\mathbb N$,
$$\mid \hat f(s)\mid = \mid \int_{2{\mathbb T}}
f(\zeta)\zeta^{-s}dm(\zeta)\mid \leq 2^{-s}
\int_{2{\mathbb T}}\mid f(\zeta)\mid dm(\zeta).$$
\end{Lem}

\begin{Lem}\label{Obs5}
For all $n\geq 0$ we have
$$\lambda F_{n+1}(\zeta)=G_{n+1,0}(\zeta)+F_{n}(\zeta^k).$$
\end{Lem}


\begin{Lem} \label{Formule2} For every $\nu \geq 1$,
$$\widehat{F_n}(\nu k^{n-1}) = \lambda^{-n}A_\nu (1-exp(-2\pi i\nu /k)).$$
In particular, 
$$\mid\widehat{F_n}(\nu k^{n-1})\mid \geq {\vert \lambda \vert}^{-n} 
\frac{\vert A_\nu \vert}{k}.$$
\end{Lem}

{\it Proof of Theorem \thesection.\ref{bori}:}

 Since $\Gamma$ acts properly 
discontinuously, the images $g_{1,n}(K)$ of the compact $K$ 
tend to infinity. Since the first component of $g_{1,n}(w,z)$
 converges, this implies that
the sequence $(\mid F_n\mid )$ converges uniformly on $K$ to $+\infty$.  \\
Since $H$ is $\mathcal C^1$ on $K$, there is a 
constant $c>0$ independent of $n$ for which

$$\beta_n:= \sup_{\zeta\in K}\mid G_{n,0}(\zeta)\mid = \sup_{\zeta\in
K}\mid H(\zeta)-H(\zeta exp(-2\pi ik^{-n})) \mid \leq ck^{-n}.$$
Fix a positive integer $N$ such that

$${\rm for\; all\;} \zeta\in
K,\;n\geq N,\mid F_n(\zeta)\mid \geq 2 \leqno{(\ast)}$$

and
$$\lambda^{-1} \sum_{n>N}\beta_n \leq 1.\leqno{(\ast\ast)}$$

Set $W:=wind(F_N)$ and  
$F(\zeta):=F_N(\zeta^k)$. Using $(\ast)$, we get $wind(F)=kW$.
We recall that
$\lambda F_{N+1}(\zeta)=F(\zeta)+G_{N+1,0}(\zeta)$.
Combine now Rouch\'e's Theorem
applied to $F$ and $G_{N+1,0}$ with 
the inequalities $(\ast)$ and 
$(\ast\ast)$ to obtain $wind(F_{N+1})=kW.$
By induction, one shows that for $p\in\mathbb N$,
$$wind(F_{N+p})=k^pW.$$
The lemmas \thesection.\ref{Obs1} and \thesection.\ref{Obs5}, imply
$$\begin{array}{lcl}
\int_{3{\mathbb T}}\ln \mid F_N(\zeta )\mid 
dm(\zeta )& =& \int_{3^k{\mathbb T}}\ln \mid F_N (\zeta )\mid 
dm(\zeta ) - (k-1) W \ln 3\\
&&\\
&=&\int_{3{\mathbb T}}\ln \mid F_N(\zeta^k )\mid dm(\zeta )- (k-1) W \ln 3\\
&&\\
&=&
\int_{3{\mathbb T}}\ln \mid \lambda  F_{N+1}
(\zeta)-G_{N+1,0}(\zeta) \mid  dm(\zeta )- (k-1) W \ln 3.
\end{array}$$

Applying  
lemma \thesection.\ref{Obs2} and the inequalities 
$(\ast)$ and $(\ast\ast)$, one gets
$$\begin{array}{lcl}
\int_{3{\mathbb T}}\ln \mid F_N(\zeta )\mid dm(\zeta )& 
\geq & \int_{3{\mathbb T}}\bigl
( \ln \mid \lambda F_{N+1}(\zeta )\mid -\mid G_{N+1,0}(\zeta)\mid\bigr) 
dm(\zeta )
- (k-1) W \ln 3\\
&&\\
&\geq & \int_{3{\mathbb T}} 
\ln \mid \lambda F_{N+1}(\zeta )\mid dm(\zeta )- \beta_{N+1}- (k-1) W \ln 3,
\end{array}$$
and furthermore
$$\int_{3{\mathbb T}} 
\ln \mid \lambda F_{N+1}(\zeta )\mid dm(\zeta ) 
\leq \int_{3{\mathbb T}}
\ln \mid F_N(\zeta )\mid dm(\zeta ) +\beta_{N+1}+ (k-1) W \ln 3.$$
Induction gives the inequalities

$$\begin{array}{lcl}
\int_{3{\mathbb T}} \ln \mid F_{N+p}(\zeta )\mid dm(\zeta )&\leq &
\int_{3{\mathbb T}}\ln \mid F_N(\zeta )\mid dm(\zeta ) +\displaystyle 
\sum_{N<s \leq N+p}\beta_s \\
&&+(k-1)W\ln 3 \displaystyle \sum_{0\leq s<p}k^s - p\, {\rm ln}\vert \lambda \vert \\
&&\\
&\leq & C + Wk^p\ln 3- p\, {\rm ln}\vert \lambda \vert,
\end{array}\leqno{(\dag\dag)}$$
for a certain constant $C>0$ independent of  $p\in\mathbb N$.

Lemma \thesection.\ref{Obs3} and 
the inequalities $(\ast)$ and $(\dag\dag)$ 
give the existence of a constant $C_1>0$ independent of $p$ 
such that
$$\ln \int_{2\mathbb T}\mid F_{N+p}(\zeta)\mid dm(\zeta) \leq C_1 + C_1k^p.$$

If $H$ {\it is not a polynomial}, 
one can find an integer $\nu $ verifying
$$\nu > C_1k^{1-N}/ \ln 2, \quad {\rm and}\quad A_\nu \neq 0.$$

But now, on the one hand Lemma \thesection.\ref{Obs4} gives
$$\begin{array}{lcl}
\ln \mid \widehat{F_{N+p}}(\nu k^{N+p-1})\mid & 
\leq & \ln \int_{2\mathbb T}\mid F_{N+p}(\zeta)\mid dm(\zeta) - 
\nu k^{N+p-1}\ln 2\\
&&\\
&\leq & C_1 + C_1k^p - \nu k^{N+p-1}\ln 2 \\
&&\\
&=&C_1 - (\nu  k^{N-1}\ln 2-C_1)k^p.
\end{array}$$

On the other hand, by Lemma \thesection.\ref{Formule2},
$$\mid\widehat{F_{N+p}}(\nu k^{N+p-1})\mid \geq {\vert \lambda \vert}^{-N-p} 
\frac{\vert A_\nu \vert}{k} >0,$$
and therefore $ \ln \vert \widehat{F_{N+p}}(\nu k^{N+p-1}) \vert \geq C_2 - C_3 p$
for some constants $C_2$ and $C_3$ independent of $p$.
This is a contradiction. \qed


\section{The contracting germ}
\setcounter{Lem}{0}

We have seen that the action of the group $\pi_{1}(S \setminus D)$
on $\mathbb H_{l} \times \mathbb C$ is generated
by the two automorphisms

$$\left\{
\begin{array}{lcl}
g_{\gamma}(w,z)&=&(w+ 2 \pi i, z)\\
g(w,z)&=&(kw,\lambda z+H(e^{-w})).
\end{array}\right.
$$

where $H(\zeta)= \sum_{m=1}^s A_{m} \zeta^m$ is
a polynomial in normal form, i.e. $A_{m}= 0$
for all $m>0$ with $k \vert m$ and $A_{s} \neq 0$.

Let $l:= [s/k]+1$.
We will conjugate our group by
$$\phi(w,z)= (w,z+\lambda^{-1} \sum_{m=1}^{l-1} A_{m} e^{-mw}).$$

\noindent This has no effect on $g_{\gamma}$, but

$$ \phi \circ g \circ \phi^{-1}(w,z) = (kw,\lambda z + Q(e^{-w})),$$
where $Q(\zeta) = H(\zeta) - \sum_{m=1}^{l-1} A_{m} \zeta^m
+\lambda^{-1} \sum_{m=1}^{l-1} A_{m} \zeta^{mk} $  is a polynomial of degree
$s$ with $$\zeta^{min(k,l)} \vert Q(\zeta).$$

Iterating this procedure if neccessary, we end up with a polynomial
$Q$ of degree $s$ such that $\zeta^l \vert Q(\zeta)$.
Let $Q(\zeta) := \sum_{m=l}^s b_{m}\zeta^m$ and
$d:= {\rm g.c.d.} \{k,\ m \mid b_{m} \neq 0\}$.

We conjugate now with $\phi(w,z) = (dw,z)$:

$$ \phi \circ g_{\gamma} \circ \phi^{-1}(w,z) = (w + 2\pi i d, z)$$

$$ \phi \circ g \circ \phi^{-1}(w,z) = (kw ,\lambda z + Q(e^{-w/d})).$$

One verifies directly that the group generated is the same
as the group $G'$ generated by

$$ (w,z) \mapsto (w+ 2\pi i, z)$$

$$(w,z) \mapsto (kw,\lambda z + Q(e^{-w/d})).$$

Let now $$l':=[s/kd] +1$$ and $$Q'(\zeta) := \sum_{m=l}^s
b_{m } \zeta^{m/d}.$$

Using the inequality 
$d[x/d]<[x]+1, d \in \mathbb N^{\star}, x \in \mathbb R$
and the fact that the indices 
of the non-vanishing coefficients of $Q'$ are divisible by
$d$,
one verifies easily that $\zeta^{l'} \vert Q'(\zeta)$.
We conjugate now with $ (w,z) \mapsto (w, e^{l'w}z)$
and the generators of $G'$ become

 $$ (w,z) \mapsto (w+ 2\pi i, z)$$

$$(w,z) \mapsto (kw,\lambda e^{l'(k-1)w} z + P(e^{w})),$$

where $P$ is the polynomial defined by

$$P(\xi):= \xi^{l'k} Q'(\xi^{-1}).$$

Remark that ${\rm deg} P \leq l'(k-1)$ and that $P(0) = 0$.
Let $$ P(\xi) = \sum_{m=1}^{l'(k-1)} c_{m} \xi^m .$$
We have ${\rm g.c.d.} \{k, m \mid c_{m}\neq 0 \} =1$.

This relation implies that the contracting germ

$$f: \Delta^{*} \times \mathbb C \rightarrow \Delta^{*} \times \mathbb
C,$$

$$f(\xi,z):= (\xi^k,\lambda \xi^{l'(k-1)}z + P(\xi))$$

is locally injective around $(0,0)$:
If $f(\xi_{1},z_{1})=f(\xi_{0},z_{0})$, then
$\xi_{0}^k = \xi_{1}^k$ and $\lambda \xi_{0}^{l'(k-1)}z _{0}+ P(\xi_{0}) =
\lambda \xi_{1}^{l'(k-1)}z _{1}+ P(\xi_{1})$.
Put $ \epsilon:= \xi_{1}/\xi_{0}$. One  has $\epsilon^k =1$ and

$$ z_{1}= \epsilon^{l'} \bigl[z_{0} +\lambda^{-1} \xi_{0}^{-l'(k-1)} 
\sum_{m=1}^{l'(k-1)}
c_{m} \xi_{0}^m (1- \epsilon^m)\bigr]. $$

\noindent If $\epsilon = 1$, one has $z_1 = z_0$.
Otherwise take $m_0$ the smallest index such that
$c_{m_0} \neq 0$ and $\epsilon^{m_0} \neq 1$.
The existence of such an index is ensured by
the relation ${\rm g.c.d.} \{k, m \mid c_{m}\neq 0 \} =1$. 
We write now

$$ z_{1}= \epsilon^{l'} \bigl[z_{0} +\lambda^{-1} \xi_{0}^{-l'(k-1)+m_0} 
\bigl(c_{m_0}(1- \epsilon^{m_0})
+\sum_{m=m_0 +1}^{l'(k-1)}
c_{m} \xi_{0}^{m-m_0} (1- \epsilon^m)\bigr)\bigr] $$
and we see that for $z_0$ and $\xi_0$ sufficiently small, 
$z_1$ stays away from $0$.
The local injectivity follows now directly.\\
By Proposition 1.2.8 of
\cite{Fav}, we get that $f$ is a defining germ 
for a minimal GSS surface $S'$
whose maximal divisor we denote by $D'$.
One can verify that the quotient of $\mathbb H_{l} \times \mathbb C$
by the action of $G'$ is the same as the one of $\Delta^{*}\times
\mathbb C$ by the equivalence relation
$u_{1} \sim u_{2} : { \textrm{ there exist }} n_{1},n_{2} \in \mathbb N$
such that $f^{\circ n_{1}}(u_{1}) = f^{\circ n_{2}}(u_{2}).$

It follows by construction that $S \setminus D$ and $S' \setminus D'$
are isomorphic. Since the intersection matrices of $D$ and of $D'$
are negative-definite and neither $D$ nor $D'$ contain
exceptional curves of the first kind, 
this isomorphism is extendable to an isomorphism of 
$S$ onto $S'$. This ends the proof of the Main  Theorem. $\qed$


\begin{thebibliography}{123}

\bibitem{Bog} {\sc Bogomolov F.A.}:
Classification of surfaces of class $VII_0$ with $b\sb{2}=0$,
{\em Izv. Akad. Nauk SSSR Ser. Mat. 40 (1976), no. 2, 273--288} 
\bibitem{D1} {\sc Dloussky G.}: Structure des surfaces de
Kato, {\em M\'emoire de  la S.M.F 112.$\rm n^{\circ}14$ (1984).} 
\bibitem{D2} {\sc Dloussky G.}: Une construction \'el\'ementaire des
surfaces d'Inoue-Hirzebruch. {\em Math. Ann. 280, (1988), 663-682.} 
\bibitem{D3} {\sc Dloussky G.}: 
Compact complex surfaces with Betti numbers $b_1=1$, $b_2>0$ and finite quotients,
{\em Proceedings of the International Congress on Differential
Geometry in memory of Alfred Gray. Bilbao, September 2000, 7 pages, submitted. }
\bibitem{DO} {\sc Dloussky G., Oeljeklaus K.}: Vector fields and foliations
on surfaces of class VII${}_{0}$, {\em Ann. Inst. Fourier 49, (1999), 1503-1545}
\bibitem{DOT1} {\sc Dloussky G., Oeljeklaus K., Toma M.}: 
Surfaces de la classe VII${}_{0}$
admettant un champ de vecteurs, {\em Comm. Math. Helv. 75, (2000), 255-270}
\bibitem{DOT2} {\sc Dloussky G., Oeljeklaus K., Toma M.}: 
Surfaces de la classe VII${}_{0}$
admettant un champ de vecteurs, II, to appear in Comm. Math. Helv.,
preprint www.mathematik.uni-osnabrueck.de/preprints/shadow/calg0009.rdf.html.
\bibitem{E} {\sc Enoki I.}: Surfaces of class VII$_0$ with curves, 
{\em T\^ohoku Math. J. 33,
(1981), 453-492.} 
\bibitem{Fav} {\sc Favre, Ch.}: Classification of $2$-dimensional
contracting rigid germs,   
{\em Jour. Math. Pures Appl. 79, (2000), 475-514}

\bibitem{Hau} {\sc Hausen J.}: Zur Klassifikation glatter kompakter $\mathbb 
C^{\star}$-Fl\"achen. {\em Math. Ann. 301 (1995), 763-769.}
\bibitem{In1} {\sc Inoue M.}: On surfaces of class VII$_0$, 
{\em Invent. Math. 24 (1974), 269-310.} 
\bibitem{KA} {\sc Kato Ma.}: Compact complex manifolds containing ``global spherical
shells'', {\em Proceedings of the Int. Symp. Alg. Geometry,  Kyoto 1977. ed. Nagata M., Kinokuniya Book Store,
Tokyo 1978.}
\bibitem{Kod}{\sc Kodaira K.} On the structure of compact complex 
analytic 
surfaces I, II,  {\em Am. J. of Math.
vol.86, 751-798 (1964); vol.88, 682-721 (1966)}.
\bibitem{LYZ}{\sc Li J., Yau S.-T., Zheng F.} On projectively flat Hermitian manifolds, 
{\em Comm.
Anal. Geom. 2 (1994), no. 1, 103--109}
\bibitem{N1} {\sc Nakamura I.}: On surfaces of class $\rm VII_0$ with
curves, {\em Invent. Math. 78,(1984), 393-443.} 
\bibitem{N2} {\sc Nakamura I.}: On surfaces of class $\rm VII_0$ with
curves II, {\em Toh\^oku Math. Jour. 42, (1990), 475-516.}
\bibitem{SS} {\sc  van Straten D., Steenbrink J.}:
Extendability of holomorphic differential forms 
near isolated hypersurface singularities,
{\em Abh. Math. Sem. Univ. Hamburg 55 (1985), 97--110.}
\bibitem{T} {\sc Teleman A.}: Projectively flat surfaces and Bogomolov's 
theorem on class $VII\sb 0$ surfaces,
{\em Int. J. Math. 5, No.2, 253-264 (1994).}
\bibitem{Zaf} {\sc Zaffran D.}: Serre problem and Inoue-Hirzebruch surfaces,
{\em Math. Ann. 319, 395-420 (2001).}


\end{thebibliography}
\end{document}